\documentclass[letterpaper, 10 pt, conference]{ieeeconf}  

\IEEEoverridecommandlockouts                              

\usepackage{graphics} 
\usepackage{epsfig} 
\usepackage{times} 
\usepackage{amsmath} 
\usepackage{amssymb}  

\usepackage{algorithm}
\usepackage{algpseudocode}
\usepackage{color}
\usepackage{hyperref}

\usepackage{subcaption}
\usepackage[ 
    backend=biber,
    style=ieee,
    sorting=none]{biblatex} 
\let\labelindent\relax
\usepackage{enumitem}
\setlist[enumerate]{
                    itemindent=1.5em,
                    leftmargin=0.5em, 
                    }

\newcommand{\set}[1]{\left\lbrace #1 \right\rbrace}

\newtheorem{theorem}{Theorem}
\newtheorem{corollary}[theorem]{Corollary}
\newtheorem{lemma}[theorem]{Lemma}

\newtheorem{definition}{Definition}

\newcounter{assumption}

\newtheorem{assumption}{Assumption}

\usepackage{cleveref}
\crefname{assumption}{assumption}{assumptions}
\Crefname{assumption}{Assumption}{Assumptions}

\crefname{theorem}{Theorem}{Theorems}
\Crefname{theorem}{Theorem}{Theorems}
\crefname{lemma}{Lemma}{Lemmas}
\Crefname{lemma}{Lemma}{Lemmas}
\crefname{corollary}{Corollary}{Corollaries}
\Crefname{corollary}{Corollaries}{Corollaries}
\crefname{proposition}{Proposition}{Propositions}
\Crefname{proposition}{Propositions}{Propositions}

\title{\LARGE \bf
Several Performance Bounds on Decentralized Online Optimization are Highly Conservative and Potentially Misleading
}

\author{Erwan Meunier and Julien M. Hendrickx$^\dagger$
\thanks{$^\dagger$E. Meunier and J. M. Hendrickx are with the ICTEAM institute, UCLouvain, Louvain-la-Neuve, Belgium. E. Meunier is supported by the French Community of Belgium through a FRIA fellowship (F.R.S.-FNRS).
J. M. Hendrickx is supported by the ``SIDDARTA'' Concerted Research Action (ARC) granted by the UCLouvain. 
Email adresses: \texttt{erwan.meunier@uclouvain.be}, \texttt{julien.hendrickx@uclouvain.be}
}}

\begin{document}

\maketitle
\thispagestyle{empty}
\pagestyle{empty}


\begin{abstract}
We analyze Decentralized Online Optimization algorithms using the Performance Estimation Problem approach which allows, to automatically compute exact worst-case performance of optimization algorithms. Our analysis shows that several available performance guarantees are very conservative, sometimes by multiple orders of magnitude, and can lead to misguided choices of algorithm. Moreover, at least in terms of worst-case performance, some algorithms appear not to benefit from inter-agent communications for a significant period of time.
We show how to improve classical methods by tuning their step-sizes, and find that we can save up to $20\%$ on their actual worst-case performance regret.
\end{abstract}
\vspace{0.5em}

\section{INTRODUCTION}
\subsection{Decentralized Online Optimization}
In recent years, there has been considerable research in online distributed multi-agent optimization and its online counterpart, the Decentralized (a.k.a Distributed) Online multi-agent Optimization (DOO) framework.
As surveyed in~\cite{surveyDOO}, the DOO setting has proven its usefulness in many applications such as medical diagnosis, distributed target tracking in 2-D plane, multi-class classification or again robot formation control. 

In decentralized multi-agent optimization, each agent $i \in V = \left[n\right]$ holds a private objective function $f_i\: : \mathcal{K} \subseteq \: \mathbb{R}^d \rightarrow \mathbb{R}$ and is aimed at cooperating with other agents so as to minimize the total cost incurred in the network. That is, computing
\begin{equation}
   \min_{x \in \mathcal{K}}  \sum_{i=1}^n f_{i}(x). \label{equation:1}  
\end{equation}
For instance, the cooperation between agents occurs when an agent $i$ communicates its estimate $x_{i,t}$ to its neighbors $\mathcal{N}(i)$, receives in its turn estimates $x_{j,t}, \, \forall j \in \mathcal{N}(i)$ and averages them w.r.t. weights given by a \textit{communication matrix} $P$:
\begin{equation}\label{example_averaging}
    \overline{x}_{i,t+1} = \sum_{j \in \mathcal{N}(i)} P_{ij} x_{j,t}.
\end{equation}
This update is an example of a \textit{consensus step}, usually coming along with an \textit{optimization step}, which can take the form of a well-known \textit{descent step}  $x_{i,t+1} = \overline{x}_{i,t} - \eta_t \nabla f_{i,t}(x_{i,t})$.

There exist many alternatives to this scheme, as in~\Cref{alg:DOCG}, where instead of averaging the information from the primal space as in~\Cref{example_averaging}, we can consider the dual space, e.g. averaging gradients $\nabla f_{i,t}(x_{i,t})$ to use more advanced algorithms. 

The Decentralized Optimization framework has been extended to the Online Optimization Setting \cite{HazanOO}, leading to the DOO framework. In this respect, each agent receives a new loss function $f_{i,t}$ at each time-step $t$ thus making it able to capture a broader variety of applications, where the problem cannot be fully known in hindsight.

The most widely used performance measure for DOO methods is the so-called \textit{Individual Static Regret} (ISR) metric which for a given agent $j$ compares the total incurred cost of choosing estimates $x_t$ (including by other agents) to the total offline cost if the decision is the minimizer $x^*$, i.e. the price paid if functions are known in advance after receiving $T$ functions
\begin{equation}\label{individual_regret}
    \textbf{Reg}_j(T) = \sum_{t=1}^T \sum_{i=1}^n f_{i,t}(x_{j,t}) - \min_{x^* \in \mathcal{K}}\sum_{t=1}^T \sum_{i=1}^n f_{i,t}(x^*). \tag{ISR}
\end{equation}
An online algorithm (centralized or not) is said to be \textit{efficient} if it incurs a sub-linear regret, that is $\lim_{t \rightarrow \infty} \textbf{Reg}_j(T)/T=0$. Other performance metrics such as the \textit{Adaptative Regret}, \textit{Dynamic Regret} and \textit{Competitive Ratio} are alternatives to the ISR. Here, we retain ISR as the main performance metric as it is the most common~\cite{book_DOO, surveyDOO}.

While performance guarantees (or worst-case bounds) have been derived for many classes of algorithms, the difference between these bounds and numerical experiments suggest they could be very conservative. Moreover, an important degree of conservatism was also observed for many bounds and performance metrics in classical decentralized optimization~\cite{sebastien_paper}. High degrees of conservatism make comparing algorithms based on their efficiency hazardous, potentially hindering progress, and prevent using the bounds for optimizing an algorithm with respect to its parameters such as step-sizes.

\subsection{Performance Estimation Problem}
The recent \textit{Performance Estimation Problem} (PEP) framework~\cite{taylor2017smooth} initially designed for automatically computing worst-case performance of continuous optimization algorithms in the offline centralized case is the key-enabler of this paper and has proven very useful on a wide range of algorithms and settings (e.g. \cite{bousselmi2024interpolation, dragomir2022optimal}). Recently, Computer-Aided methods for analyzing Online Optimization Algorithms have seen a growing interest whether by using PEP on Online Frank-Wolfe~\cite{weibel2025optimized} or the Integral Quadratic Constraints framework on Online Gradient Descent~\cite{jakob2025online}.

In practice, suppose one wants to know the worst-case performance of the so-called Gradient Descent (GD) for functions lying in the set of convex and $L$-Lipschitz functions, denoted $\mathcal{F}_{L}$. This boils down to finding the function $h \in \mathcal{F}_L$ maximizing $\mathcal{P}(h,x_1,x^*) := h(x_T) - h(x^*)$ where the sequence of estimates $x_1, \ldots, x_T$ are obtained from GD:
\begin{equation}\label{update_step_GD}
    x_{t+1} = x_t - \alpha \nabla h(x_t), \quad \forall t \in \left[T-1\right]
\end{equation}
with the (classical) \textit{initial condition} $\|x_1 - x^* \|\leq D$ where $x^* = \arg \min h$ . The PEP framework tackles this challenge by formulating a $\texttt{PEP-Model}$ through the following three steps:
\\
    \textbf{Step 1: }   \textit{\textbf{Discretizing}} the function $h$ and its gradient in all estimates so as to obtain a set of decision variables $\mathcal{D}:=\set{x_t,h_t,g_t}_{t \in \left[T\right]\cup\set{*}}$ where $h_t$ represents $h(x_t)$ and $g_t$ is $ \nabla h(x_t)$. At this step, the vectors $x_t$, $h_t$ and $g_t$ are independent from each other and not \textit{a priori} consistent with an actual function.
    \\
    \textbf{Step 2:}
    Imposing the so-called \textit{\textbf{interpolation constraints}} which are necessary and sufficient conditions bending decision variables in $\mathcal{D}$, to be consistent with an interpolating function $h$ lying in a given class (e.g Lipschitz, Smooth, Strongly-Convex, etc) w.r.t. some specified parameters. 
    In our case, $h$ is imposed to be in $\mathcal{F}_L$ via the following constraints (cf. Theorem 3.4 of \cite{taylor2017smooth}) on $\mathcal{D}$; \begin{equation}\label{PEP_interpol_conditions}
    \begin{aligned}
        h_{i} &\geq h_j + \langle g_j, x_i - x_j \rangle \\
        \| g_j \| & \leq L 
    \end{aligned}
    \quad , \forall i, j \in \left[T\right]\cup\set{*},
    \end{equation}
    which allows keeping the equivalence with the initial problem of maximizing $\mathcal{P}(h,x_1,x^*) := h(x_T) - h(x^*)$.
    \\
\textbf{Step 3:}
    In addition to interpolation constraints, we further constrain estimates so as to encode; the \textit{performance metric} to be maximized
    \begin{equation}\label{PEP_objective}
    \sup_{\set{x_t,h_t,g_t}_{t\in\left[T\right]\cup\set{*}}}  h_T - h^*,
    \end{equation}
    \textit{update steps} representing~\Cref{update_step_GD}
    \begin{equation}\label{PEP_update_step}
        x_{t+1} = x_{1} - \alpha \sum_{k=1}^t  g_k, \quad \forall t \in \left[T-1\right]
    \end{equation}
    as well as \textit{initial} and \textit{optimality} conditions: 
    \begin{equation}\label{PEP_init_opt_conditions}
        \|x_1 - x^*\| \leq D  \text{ and } g^* = 0. 
    \end{equation}
Finally, the \texttt{PEP-model} is obtained by wrapping up~\Cref{PEP_interpol_conditions,PEP_objective,PEP_update_step,PEP_init_opt_conditions}.

It is worth mentioning that, under common assumptions (e.g. $\alpha$ being fixed), \texttt{PEP-model} can be reformulated into a Semidefinite Program (SDP).

In~\cite{sebastien_paper} PEP has been extended to the Decentralized Offline Multi-agent Optimization setting, which shown dramatic improvements regarding the bounds, sometimes by orders of magnitude, especially with respect to the network topology. Furthermore, \cite{colla2024exploiting} showed that, up to a rescaling by the number of agents, many worst-case performance metrics for standard algorithms were independent of the number of agents.
Taking a step back, PEP is often used to choose the best method w.r.t. best (near)-tight worst-case bounds and to improve / design the most promising method by acting on its parameters (e.g. step-sizes). 
\subsection{Contributions}
In light of the recent developments 
regarding both DOO algorithms and the PEP framework our contributions unfold in 5 points:
\begin{enumerate}[wide,labelwidth=!, labelindent=0pt] 
    \item Formulating DOO algorithms in PEP, which is directly allowed by the current decentralized PEP framework.
    \item Providing (near)-tight worst-case bounds for DOO algorithms within the PEP framework. Notably, we show that several analytical bounds are very conservative and can mislead a user into picking the ``best method''. \label{contribution:1}
    \item Demonstrating that, in some cases, a projection-free algorithm derives little benefit from communication among agents in the worst-case and for several steps. 
    \item Showing how we can select methods in the light of b).
    \item Improving DOO methods by finding better step-sizes using a black-box optimization algorithm. Regarding c), we find step-sizes making projection-free algorithms benefiting more from the communication.
\end{enumerate}

\section{Problem Setting}
\subsection{Problem and notations}
\subsubsection{Distributed Online Optimization}
In this article, we focus on finding a set of initial points and a sequence of functions maximizing the performance of DOO algorithms aimed at minimizing (\ref{individual_regret}) over a static network $\mathcal{G}=(V,E)$, where $V=\set{1,\ldots,n}$ is the set of agents
and $E=\set{(i,j) \mid P_{ij}\neq0, i, j \in V}$ is the set of activated edges, where $P \in \mathbb{R}^{n \times n}$ is the weights matrix that represents the weighted ``belief'' between agents.

The goal of every node $i \in V$ in the network is to generate a sequence of estimates $\set{x_{i,t}}_{t=1}^T \in \mathcal{K}$, from an algorithm $\mathcal{A}$, aimed at minimizing (\ref{individual_regret}) where $T\in\mathbb{N}^*$ is the number of functions received at the end of the day (also the number of iterations for DOO methods) and $f_{i,t}$ are defined over a shared closed convex set $\mathcal{K}$. We denote by $\nabla f_{i,t}(x)$ the (sub)gradient of $f_{i,t}$ evaluated in $x \in \mathcal{K}$, and as a shorthand notation $g_{i,t} := \nabla f_{i,t}(x_{i,t})$.
\subsubsection{PEP}
Due to some limitations over the representability of non-negative matrices within PEP discussed in~\cite{sebastien_paper}, a slightly more general class of matrices is introduced.
\begin{definition}[Generalized doubly stochastic matrix~\cite{chiang2005generalized}] \label{def_GDSM}
    A matrix $P \in \mathbb{R}^{N \times N}$ is \textit{generalized doubly stochastic} if its rows and columns sum to one, i.e., if 
    \begin{equation*}
        \sum_{i=1}^n P_{ij} = 1, \quad \sum_{i=1}^n P_{ji} =1, \quad \forall j \in V.
    \end{equation*}
\end{definition}
The resulting worst matrix of the PEP may thus have negative elements. However, the provided worst-case guarantees hold for non-negative matrices. 
More, eigenvalues of $P$ are sorted as follows $-1 \leq \lambda_p \leq \ldots \leq \lambda_2 \leq \lambda_1=1$, where $\lambda_2$ is the second largest eigenvalue of $P$ and $1\leq p \leq n$. $\lambda_2$ can be interpreted as the ``network disconnectedness''. 

We denote by $\mathcal{W}(\lambda^-, \lambda^+)$ the set of real, symmetric and generalized doubly stochastic $n \times n$ matrices having their eigenvalues between $\lambda^-$ and $\lambda^+$ except for $\lambda_1=1$; $
    \lambda^- \leq \lambda_N \leq \cdots \leq \lambda_2 \leq \lambda^+$ where $  \lambda^-, \lambda^+ \in (-1,1)
$. 

In what follows, by abuse of notation
$\| \cdot\|$ denotes the $\ell_2$-norm.
$\mathcal{F}$ is an arbitrary class of functions and $\mathcal{A}$ is an arbitrary DOO algorithm compatible with previous assumptions.

\subsection{Shared Assumptions}
We now detail the assumptions shared among algorithms in this paper, which are very common among works in DOO~\cite{surveyDOO}. Assumptions specific to each algorithm are detailed in dedicated sections.
\begin{assumption}[Network topology]
    The network $\mathcal{G}$ and the weight matrix $P$ satisfy the following.
\begin{enumerate}
    \item $\mathcal{G}$ is a strongly connected graph, meaning that there exists a path between every pair of agents.
    \item $P$ is a generalized doubly stochastic as in~\Cref{def_GDSM}.
\end{enumerate}
\end{assumption}

\begin{assumption}[Bounded set]\label{assumption:bounded_set}
    The decision space $\mathcal{K}$ has a bounded diameter; that is, $\forall x, y \in \mathcal{K}, \|x-y\|\leq D$. 
\end{assumption}

\begin{assumption}[Lipschitz continuous functions]\label{assumption:Lipschitz_continuous_function}
    Function $f_{i,t}$ is $L$-Lipschitz over $\mathcal{K}$, i.e., $\left| f_{i,t}(x) - f_{i,t}(y) \right| \leq L \|x-y\|$ for all $x, y \in \mathcal{K}, i \in V, t \in \left[T\right]$.
\end{assumption}


\section{Analyzed DOO algorithms}

We introduce the so-called Distributed Autonomous Online Learning (DAOL)~\cite{DAOL_distributed_autonomous_online_learning}, Distributed Online Conditional Gradient (DOCG)~\cite{DOCG_distributed_online_conditional_gradient} and Distributed Online Mirror Descent (DOMD)~\cite{DOMD_distribute_online_mirror_descent}. 
Despite not analyzed in the present, we note the existence of~\cite{Dual_averaging} which is a specific case of DOMD, i.e., our analysis for DOMD can be used to directly assess performances of~\cite{Dual_averaging}. The scale-invariance property w.r.t. $L$ and $D$ is discussed in~\Cref{section:scale-invariance}, which notably allows to compare algorithms on a fair basis.

\subsubsection{Distributed Autonomous Online Learning}
Distributed Autonomous Online Learning (DAOL), also known as Distributed Online Gradient Descent is a natural extension of the so-called \textit{Online Gradient Descent} (OGD) and \textit{Decentralized Gradient Descent }(DGD). This algorithm is the basis for more advanced algorithms and takes advantages of its simplicity to be understood and implemented. It comes with very few assumptions on $f_{i,t}$.
\begin{assumption}[Identical initialization for estimates]
    Initialize identically the first estimates $x_{i,1} \in \mathcal{K}, \forall i \in V$.
\end{assumption}
A description of DAOL is given in \Cref{alg:DAOL} and comes along with an upper-bound over ISR in~\Cref{theorem_DAOL}.

\begin{algorithm}[h]
    \caption{DAOL}
    \label{alg:DAOL}
    \begin{algorithmic}[1]
        \For{$t \in \left[T\right]$ and each learner $i \in V$}
                \State $z_{i,t+1} \gets \sum\limits_{j \in \mathcal{N}(i)}  P_{ij} x_{j,t} - \eta_t g_{i,t}$ \hfill (Consensus step)
                \State $x_{i,t+1} \gets \arg\min\limits_{w \in \mathcal{K}} \| w - z_{i,t+1} \|$ \hfill (Projection step)
        \EndFor
    \end{algorithmic}
\end{algorithm}

\begin{theorem}[Bound over ISR for DAOL \cite{DAOL_distributed_autonomous_online_learning}]\label{theorem_DAOL}
Pick parameters $\eta_t=\frac{D}{2L\sqrt{t}}$, then 
\begin{equation*}
    \mathbf{Reg}_j(T) \leq \left(5+\frac{16}{1-\lambda_2}\right)nLD\sqrt{T}.
\end{equation*}
\end{theorem}
\vspace{0.5em}
\textbf{Remark:}
We note that~\cite{DAOL_distributed_autonomous_online_learning} propose parameters which enable DAOL to handle the strongly-convex case, we however decide to focus on the non-smooth setting.
\subsubsection{Distributed Online Conditional Gradient Descent}
As in the centralized offline setting, the projection step can sometimes be costly to compute \footnote{For example, sparse optimization over the $\ell_1$-Ball.}. Whence using Frank-Wolfe conditional steps can help if the feasible set has some good properties. DOCG leverages this idea in the context of DOO at the cost of a worse dependence on time.
\begin{assumption}[Network topology for DOCG]
    The communication matrix $P$ used in DOCG is symmetric.
\end{assumption}

\begin{assumption}[Initial conditions for DOCG]
    Fix $x_{i,1} \in \mathcal{K}$ and $ z_{i,1} = 0$, for all $i \in V$.
\end{assumption}
DOCG is depicted in~\Cref{alg:DOCG} and an upper-bound on its ISR is given in~\Cref{bound_DOCG}.
\begin{algorithm}[h]
    \caption{DOCG}
    \label{alg:DOCG}
    \begin{algorithmic}[1]
        \For{$t \in \left[T\right]$ and each learner $i\in V$}
                \State \textbf{Define} $F_{i,t}(x) := \eta_i \langle z_i(t), x \rangle + \|x - x_{1,1}\|^2$
                \State $v_{i,t} = \arg\min\limits_{w \in \mathcal{K}} \{ \langle \nabla F_{i,t}(x_{i,t}), w \rangle \}$ \hfill \text{(Lin. Opt. Step)}
                \State $x_{i,t+1} = x_{i,t} + \gamma_{t,i} (v_{i,t} - x_{i,t})$ \hfill \text{(Conv. Combination)}
                \State $z_{i,t+1} = \sum\limits_{j \in \mathcal{N}(i)} P_{ij} z_{j,t} + g_{i,t}$ \hfill (Dual Consensus Step) 
        \EndFor
    \end{algorithmic}
\end{algorithm}

\begin{theorem}[Bound over ISR for DOCG~\cite{DOCG_distributed_online_conditional_gradient}]\label{bound_DOCG}
    The DOCG algorithm with parameters $\eta_i = \frac{(1-\lambda_2)D}{2(\sqrt{n}+1+(\sqrt{n}-1)\lambda_2)LT^{3/4}}$ and $\gamma_{t,i} = \frac{1}{\sqrt{t}}$ for all $i \in V$ and any $t\in\left[T\right]$ attains the following regret bound:
    \begin{equation}\label{bound_DOCG_eq}
    \begin{aligned}
        \mathbf{Reg}_j(T) \leq & \left(8n + 2\left(\frac{1+\lambda_2}{1-\lambda_2}\right)\sqrt{n}\right) LDT^{3/4}+ \beta LDT^{1/4}
    \end{aligned}
     \end{equation}
     where 
    $
         \beta  \in \left[\frac{1}{2},\frac{3}{2}+\frac{1}{8\sqrt{n}}\right]\subseteq\left[\frac{1}{2},\frac{13}{8}\right]
    $
     is a constant becoming negligible with respect to the leading term in $\mathcal{O}\left(T^{3/4}\right)$.
\end{theorem}
\textbf{Remark:}
We notice that regularization term used in line 2 of \Cref{alg:DOCG} bias the towards agent 1's first estimate. However authors of~\cite{DOCG_distributed_online_conditional_gradient} do not assume that all initial estimates are the same. In our opinion, the regularization term could be replaced by $\|x - x_{i,1} \|^2$ which does not require any previous communication of agent 1's first estimates among the network. It is further discussed in~\Cref{par:impact_of_the_communication_network}.
\subsubsection{Distributed Online Mirror Descent}

Here we study a special case\footnote{Contrary to~\cite{DOMD_distribute_online_mirror_descent}, the network is considered static and there is no approximation errors.} of the more general algorithm Distributed Mirror Descent for Online
Composite Optimization introduced in~\cite{DOMD_distribute_online_mirror_descent} under the name of Distributed Online Mirror Descent (DOMD),
akin to~\cite{book_DOO}.
DOMD generalizes DAOL by using the Bregman Divergence $V_\omega$ which is defined as
\begin{equation*}
V_{\omega}(x,y) := \omega(x) - \omega(y) - \langle  \nabla\omega (y), x -y \rangle
\end{equation*}
with respect to the Kernel function $\omega:\mathbb{R}^d \longrightarrow \mathbb{R}$ complying to the following assumptions:
\begin{assumption}
The kernel $\omega$ is $\sigma_\omega$-strongly convex i.e. $\omega(x) - \frac{\sigma_\omega}{2}\| x \|^2$ is convex.
\end{assumption}
\begin{assumption}
The kernel $\omega$ is $G$-smooth i.e. $\| \nabla \omega (x) - \nabla \omega(y) \| \leq G \|x-y\|$. 
\end{assumption}

For the rest, we do not ask for more assumptions on $V_{\omega}$. Note, that the work by~\cite{dragomir2022optimal} enables us to analyze such operators in PEP.
\begin{assumption}\label{assumption:epsilon}
     $\exists \varepsilon>0$ such that $P_{ij} \geq \varepsilon$ for all $(i,j) \in E$.  
\end{assumption}
\begin{assumption}[Initial conditions for DOMD]
     Fix the initial estimates to be $x_{i,1} = \arg\min_{x \in \mathcal{K}} \omega(x)$ for all $i \in V$.
\end{assumption}
DOMD is given in~\Cref{alg:DOMD} and the bound over the incurred ISR is given by~\Cref{theorem_DOMD}.
\\
\textbf{Remark:}
We note that, authors of~\cite{DOMD_distribute_online_mirror_descent} do not assume that the communication matrix $P$ is symmetric for DOMD. Our results, based on PEP are thus not as general as the one given in~\Cref{theorem_DOMD} but will still hold for a wide-variety of communication networks.
Also,
\Cref{assumption:epsilon} makes DOMD not directly translatable in PEP, we discuss and overcome this limitation in Lemma~\ref{lemma:epsilonlambda}.
The latter allows us to compare algorithms on a common ground for $N=2$ agents by varying $\lambda_2$.

\begin{lemma}[Upper-bound for the second largest eigenvalue]\label{lemma:epsilonlambda}
    If $P$ is of size $2$ and $P_{ij} \geq \varepsilon, \, \forall i,j \in V$ if and only if $ \lambda_2 \leq 1-2\varepsilon$.
\end{lemma}

\begin{proof}
    $P$ is a doubly stochastic matrix, then it as the form, $ P = \left(\begin{smallmatrix}
        a & 1-a \\
        1-a & a
    \end{smallmatrix}\right)$ whose eigenvalues are $\lambda_1 = 1$ and $\lambda_2 = 1-a$. Then, using the lower bound over entries we have $a \in \left[\varepsilon,1-\varepsilon\right]$ from which $\lambda_2 \in \left[2\varepsilon, 1-2 \varepsilon\right]$. 
\end{proof}

\begin{algorithm}[h]
    \caption{DOMD}
    \label{alg:DOMD}
    \begin{algorithmic}[1]
        \For{$t \in \left[T\right]$ and each learner $ i \in V$}
                \State $z_{i,t} = \arg\min\limits_{w \in \mathcal{K}} \lbrace \langle g_{i,t}, w \rangle + \frac{1}{\eta} V_{\omega}(w, x_{i,t}) \rbrace$ \hfill (Mir. Step) 
                \State $x_{i,t+1} = \sum\limits_{j \in \mathcal{N}(i)}  P_{ij} z_{j,t}$ \hfill (Consensus Step)
        \EndFor
    \end{algorithmic}
\end{algorithm}
\begin{theorem}[Bound over ISR for DOMD \cite{DOMD_distribute_online_mirror_descent}]\label{theorem_DOMD}
Pick parameter $\eta > 0$. Then for all $T\geq 1$ and any $j \in V$, there holds
    \begin{equation}\label{eq:Regret_DOMD}
        \mathbf{Reg}_j(T) \leq A_0 + \frac{A_1}{\eta}+A_2T\eta \quad \text{where}
    \end{equation}
    \vspace{-1em}
    \begin{align*}
        A_0 = 2L\phi\sum_{i=1}^n  \| x_{i,1} \|, &\hspace{2em} A_1 = \sum_{i=1}^n V_{\omega}(x^*,x_{i,1}), \\
        A_2 = \frac{nL^2(1+2\phi)}{2\sigma_{\omega}}, &\hspace{1.9em} \phi = \frac{64n^6}{\varepsilon^3 + 16 \varepsilon n^4-8(\varepsilon n)^2}.
    \end{align*}
\end{theorem}
Later on, we use $x_{i,1}$ and $x^*$ returned by PEP to compute the bound presented in~\Cref{theorem_DOMD}. When $n>2$, $\varepsilon$ is computed by taking the largest non-null value of the estimated worst-case communication matrix returned by $P$. A simpler bound can also be obtained by bounding $\left\| x_{i,1}\right\|$.
\begin{corollary}[\cite{DOMD_distribute_online_mirror_descent}]
    If $\eta = \frac{D}{L\sqrt{T}}$, then the ISR of DOMD is bounded as follows:
    \begin{equation}\label{eq:bound_DOMD}
        \mathbf{Reg}_j(T) \leq A_0 + \left(\frac{L}{D}A_1+\frac{D}{L}A_2\right) \sqrt{T} .
    \end{equation}
\end{corollary}

\subsubsection{Scale-invariant Algorithms}\label{section:scale-invariance}

We note that both DAOL and DOMD as described in their respective seminal papers are not scale-invariant. 
For instance, notice that running ~\Cref{alg:DAOL} or~\Cref{alg:DOMD}, on $\Tilde{f}(x) := 2f(x)$ and $\Tilde{L}:=2L$ does change the trajectory taken by the algorithm. 
That is why \Cref{theorem_DAOL,theorem_DOMD} as depicted in the present paper have been slightly re-written (both papers give a worst-case bound w.r.t. the step-sizes) to take the scaling into consideration. 
In practice, in both cases we have multiplied their respective step-sizes strategy by $\frac{D}{L}$. For DOMD, we directly obtain the scaled version of their bound by plugging $\eta=\frac{D}{L\sqrt{T}}$ into~\Cref{theorem_DOMD} and verifying that~\Cref{eq:bound_DOMD} is proportional to $LD$.
For DAOL, one can remark that instantiating their main bound for $D=L=1$ in~\cite{DAOL_distributed_autonomous_online_learning} and then multiplying again by $\tilde{L}\tilde{D}$ yields a linear scaled expression. 
We refer to~\cite{colla2024exploiting} for an in-depth discussion about the \textit{scale-invariance}.

\section{PEP model for DOO algorithms}
Following the same approach as in our example where we apply the PEP framework to GD (see ~\Cref{PEP_objective,PEP_interpol_conditions,PEP_update_step,PEP_init_opt_conditions}), we define the decision variables of the PEP formulation for representing iterations of DOO algorithms. In this respect, we define
 $X := \set{x_{i,t} \in \mathbb{R}^d \mid i \in V, t \in \left[T\right]}\cup \set{x^* \in \mathbb{R}^d}$ the set of all estimates, 
 $\Phi := \set{f_{i,t}^{j} \in \mathbb{R}}\cup\set{f_{i,t}^* \in \mathbb{R}}$ the set of functions values evaluated in all estimates, and by the same token, their respective gradients $\Gamma :=\set{g_{i,t}^{j}\in \mathbb{R}^d}\cup\set{g_{i,t}^*\in\mathbb{R}^d}$ where $d$ is the dimension of the worst-case instance\footnote{Formulating a PEP into a SDP removes the dependency of the worst-case instance on $d$ (c.f.  Chapter 4.2 of~\cite{taylor2017smooth}).}. More, denote $\mathcal{D}=X\cup \Phi \cup  \Gamma $ the set of all PEP decision variables, for the sake of compactness.

The high-level PEP formulation for DOO algorithms whose performance metric is the individual regret (\ref{individual_regret}), is given by~\ref{PEP_doo}.
\begin{equation}\label{PEP_doo}
\begin{aligned}
                \underset{\mathcal{D}}{\sup} &  \sum_{t=1}^T\sum_{i=1}^n f_{i,t}^j - f_{i,t}^*\\
                & x_{i,t} \text{ estimates obtained from $\mathcal{A}$ applied on $f_{i,t}$}\\
                & \text{Consensus steps w.r.t. } \mathcal{A} \text{ using } P ,\\
                & \text{Interpolation constraints on } \mathcal{D}, \\
                &\text{Initial conditions on } \mathcal{D}, \\
                &  x^*\in \arg \min_{\mathcal{K}} \sum_{t=1}^T \sum_{i=1}^n f_{i,t} .
\end{aligned} \tag{\texttt{PEP-DOO}}
\end{equation}
We refer to \cite{taylor2017convex} for a detailed way of reformulating \ref{PEP_doo}
into an SDP. 

In addition, representing the consensus step and assumptions over $P$ in PEP can be found in~\cite{sebastien_paper}. About the latter point, $P$ can be either specified explicitly (i.e. all entries are given) or implicitly, by constraining the range of the second largest eigenvalues, i.e. assuming that $\lambda_2 \in \mathcal{W}(\lambda^-, \lambda^+)$. If the \textcolor{blue}{exact} communication matrix is provided, the resulting PEP worst-case bound is \textcolor{blue}{\textbf{tight}}, otherwise when considering a \textcolor{red}{range} of $\lambda_2$ for $P$, the result will be \textcolor{red}{\textbf{near tight}}, since~\cite{sebastien_paper} used a relaxation so to represent the range of eigenvalues. However, if the matrix range is used,~\cite{sebastien_paper} showed how to recover the worst-case communication matrix up to a small-error in most cases, thus indicating that the worst-case bound is guaranteed to be tight or not.  

\section{PEP applied to DOO algorithms} 
Our results were obtained by using the PESTO
toolbox~\cite{PESTO} which is a MATLAB package aimed at automatically generating SDP formulations for PEP problems (e.g. \ref{PEP_doo} in our case).

Without loss of generality, since the analyzed algorithms are scale-invariant, a normalized version of ISR is defined as follows:
\begin{equation}
\overline{\mathbf{Reg}}_j(T) = \mathbf{Reg}_j(T)/TnLD, \tag{\text{N-ISR}}
\end{equation}
and called the \textit{Normalized Individual Static Regret} (N-ISR) so to make performances comparable.
Observe that the worst-possible case is to have an algorithm showing a regret amounting to $TnLD$, whence $\overline{\mathbf{Reg}}_j(T)\in\left[0,1\right]$.

\subsection{Set-up}\label{sub_section:set_up}

In Fig~\ref{fig:RegwrtLambda2}, we represent on top, the N-ISR bounds coming from the Literature and the bottom figure depicts the near tight PEP bounds, against $\lambda_2 \in \left[0,1\right]$.  Fig~\ref{fig:RegwrtTandN} plots the bounds for N-ISR, with respect to $T$ on the left-hand side, and with respect to $N$ on the right-hand side. As in Fig~\ref{fig:RegwrtLambda2}, since the Literature bounds are always greater than PEP bounds by at least $1$ order of magnitude, we plot the Literature bounds on the top and the PEP bound on the bottom.

For all figures, if not picked as a variable for the $x$-axis \textbf{we fix}:
$\boxed{n=2, \, D=1, \,L=1,\, G=2, \, \sigma_\omega = 1, \,\lambda_2=0.9}$ (resp. Number of agents, Diameter, Lipschitz constant, Smoothness constant for the Kernel function, Strong-convexity constant for the Kernel function used in DOMD, Second largest eigenvalue of the communication matrix). 
$T$ and $n$ are kept relatively moderate in size to ensure that numerical computations are performed within a reasonable amount of time (few minutes).
Also, by abuse of notation in this section, fixing $\lambda_2$ means that we consider $P \in \mathcal{W}(-\left|\lambda_2\right|, \left|\lambda_2\right|)$.
Furthermore, since the regret is dependent on the choice of an agent, results for $n=2$ agents cannot be extended in the same fashion as in ~\cite{colla2024exploiting}. However, in the later we will show that informative tendencies can be observed and even leveraged for improving performances.

\subsection{Comparison of PEP and Analytical bounds}

\begin{figure}[htpb]
    \centering
    \includegraphics[width=0.9\linewidth]{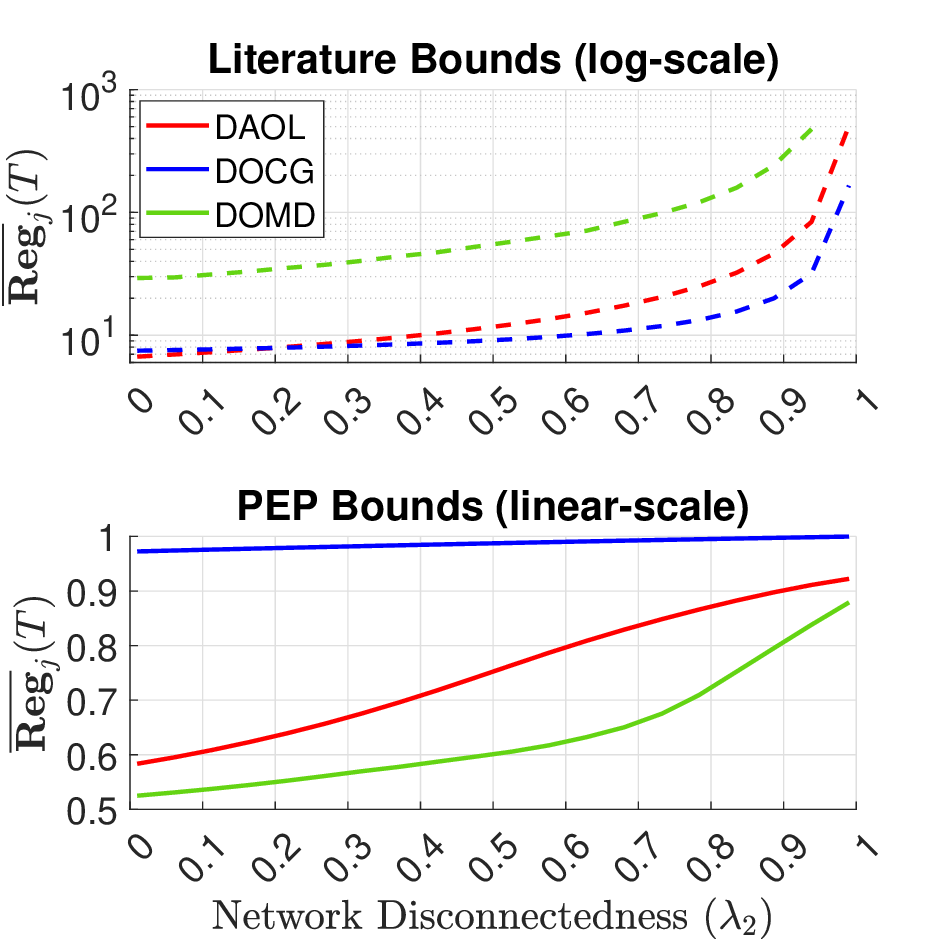}
    \caption{Comparison of analytical bounds from the Literature and near tight bounds given by PEP w.r.t. to the network topology controlled by $\lambda_2$ (with $T=10$).\\} 
    \label{fig:RegwrtLambda2}
\end{figure}

\begin{figure}[htpb]
    \centering
    \setkeys{Gin}{width=0.51\linewidth} 
\subfloat{\includegraphics{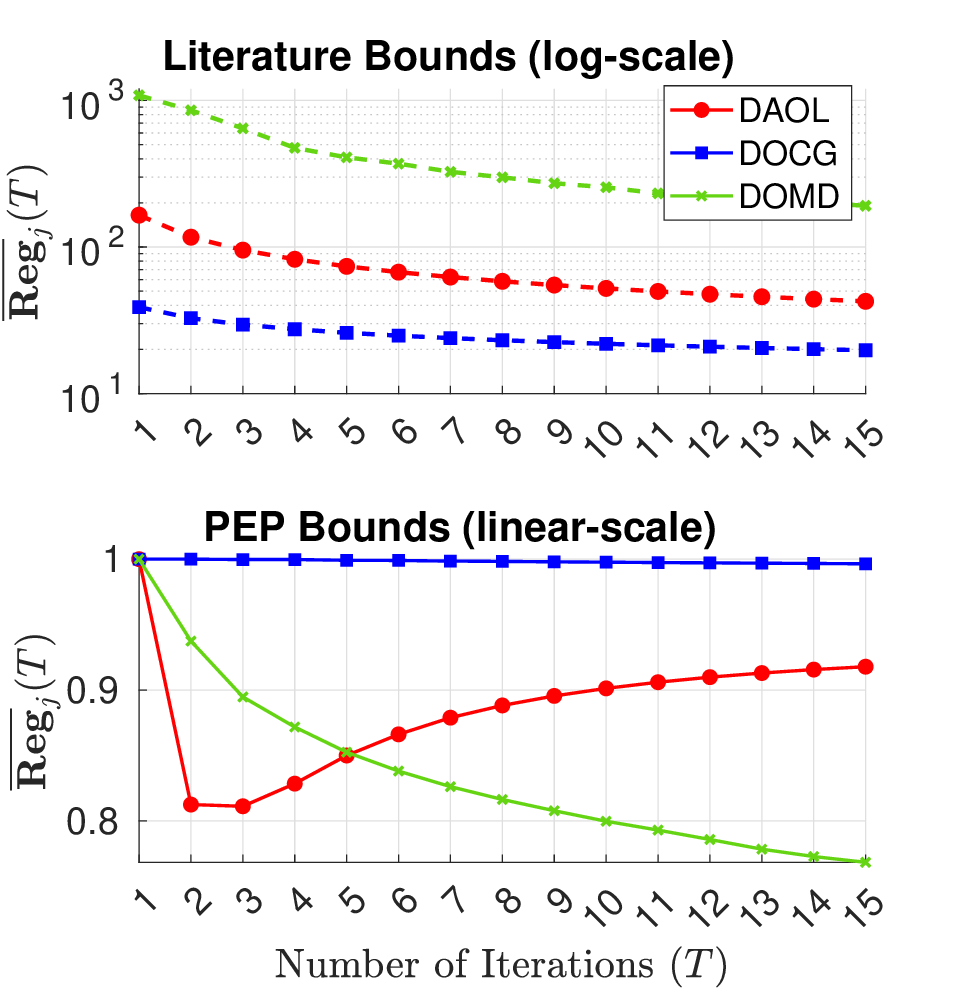} }
\subfloat{\includegraphics{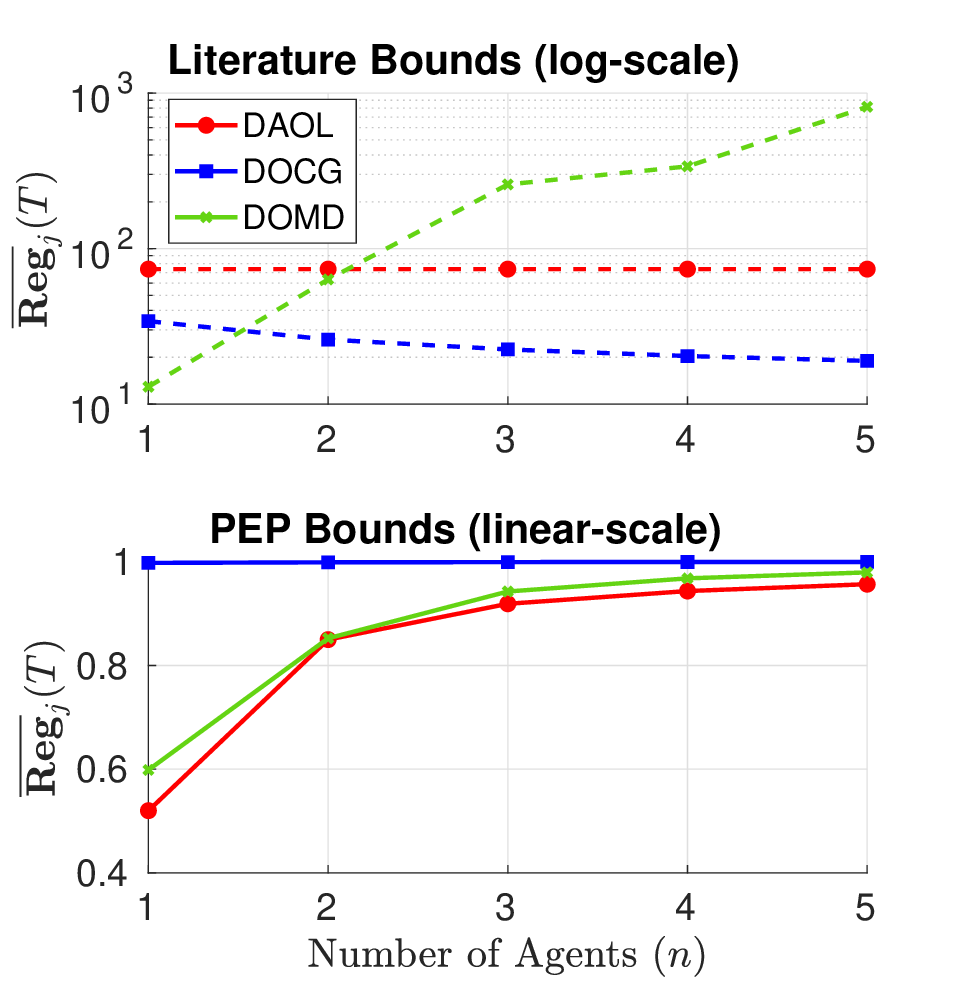} }

\caption{Comparison of analytical bounds from the Literature and near-tight bounds given by PEP against the number of iterations $T$ on the left, and the number of agents $n$ on the right (with $T=5$).}
\label{fig:RegwrtTandN}
\end{figure}
First, we notice that near tight bounds returned by PEP are significantly lower \--- by up to $2$ orders of magnitude \--- than those from the Literature for all studied algorithms and parameters taken into account ($\lambda_2$ in Fig. \ref{fig:RegwrtLambda2}, $T$ and $n$ in Fig. \ref{fig:RegwrtTandN}). In Fig.~\ref{fig:RegwrtLambda2}, we see that the N-ISR given by analytical bounds for all algorithms increase and even diverge as the connectivity of the network decreases (i.e. $\lambda_2$ goes to $1$). In other words, analytical bounds do not predict accurately the performance of algorithms on poorly connected networks\footnote{The sequence of second largest eigenvalues $\lambda_2(n)$ of a $n\times n$ grid graph (each node is connected to itself and its $4$ surrounding neighbors with weight $0.2$) for $N=2$ to $6$ are; $\lambda_2(n)\simeq0,0.4,0.56,0.68,0.77$. Poorly connected graph are thus typical.}. For instance, analytical bounds coming from the Literature suggest that, on our class of 
problems, the DOCG method is the most advisable.  

Moreover, \cite{DOCG_distributed_online_conditional_gradient} claim the superiority of DOCG over DAOL w.r.t. network topology. We see, in fact, that despite analytical upper-bounds support their claim; near-tight PEP bounds indicate that using DAOL or even more DOMD turns to be a more reasonable choice for tested parameters and w.r.t. $N$, $T$ or again $\lambda_2$.
\subsubsection{DAOL}
In~\Cref{fig:RegwrtTandN}, the regret seems to evolve sub-linearly with $T$. However, we observe on better connected graphs a decrease happening for small $T$. We suspect that $\overline{\mathbf{Reg}}_j(T)$ stabilizes and decreases after a reasonable number of iterations for poorly connected graphs.
\subsubsection{DOCG}\label{par:impact_of_the_communication_network}
\paragraph{Impact of the communication network}
Surprisingly, in Fig.~\ref{fig:RegwrtLambda2}  we see that the worst-case bound returned by PEP for DOCG indicates that the (N)-ISR, (which is very close to $1$) barely depends on $\lambda_2$. In other words, in the worst-case scenario and in the range of tested parameters, DOCG does not appear to benefit from communication among agents towards reaching a better consensus.
On this specific issue, we have performed numerous computations of PEP bounds in various settings. Our results were validated in an ad-hoc way by  checking that estimates, gradients and function values returned by PEP (i) verify assumptions (w.r.t. $L$ and $D$) as well as (ii) follow DOCG. 
This was verified for \(T \in [50]\), though it may only hold within the tested range of parameters. This does not invalidate~\Cref{bound_DOCG}, as PEP bounds remain below it, and the algorithm may produce significant results only after many iterations. In practice, for \(T \ll 3.45 \cdot 10^4\), \cref{bound_conjecture} is less than or equal to \cref{bound_DOCG_eq}, since the latter contains large constants\footnote{Even with \(\lambda_2 = 0.6\), \(L = D = 1\), and \(n = 2\), which minimize these constants.}. Although \Cref{bound_DOCG} via \cref{bound_DOCG_eq} guarantees that DOCG runs in at most \(\mathcal{O}(T^{3/4})\), while \cref{bound_conjecture} scales as \(\mathcal{O}(T)\), this is not contradictory. We hypothesize two regimes: one where our bound holds for moderate \(T\), and an asymptotic regime governed by \Cref{bound_DOCG}. Furthermore, for DOCG, the regret appears to grow linearly with both the number of agents \(n\) and the number of iterations \(T\). Our numerical experiments suggest that, for \(\lambda_2 \gtrsim 0.6\), the ISR for DOCG satisfies
\vspace{-0.3em}
\begin{equation} \label{bound_conjecture}
     \mathbf{Reg}_j(T) \propto T n L D,
\end{equation}
with step-sizes chosen as in~\Cref{theorem_DOMD}. We do not observe similar behaviour for the other two methods analyzed, which both benefit from inter-agent communication.

\paragraph{On the regularization}
Moreover, in the continuation of our remark about line 2 of~\Cref{alg:DOCG}, changing the regularization from $\|x - x_{\textcolor{red}{1},1}\|^2$  to what we think would be a better choice (since it requires weaker and more realistic assumptions), i.e., $\| x - x_{\textcolor{blue}{i},1}\|^2$ for $i \in V$ has left PEP bounds unchanged. 
\subsubsection{DOMD}
\begin{figure}[htpb]
    \centering
    \includegraphics[width=0.9\linewidth]{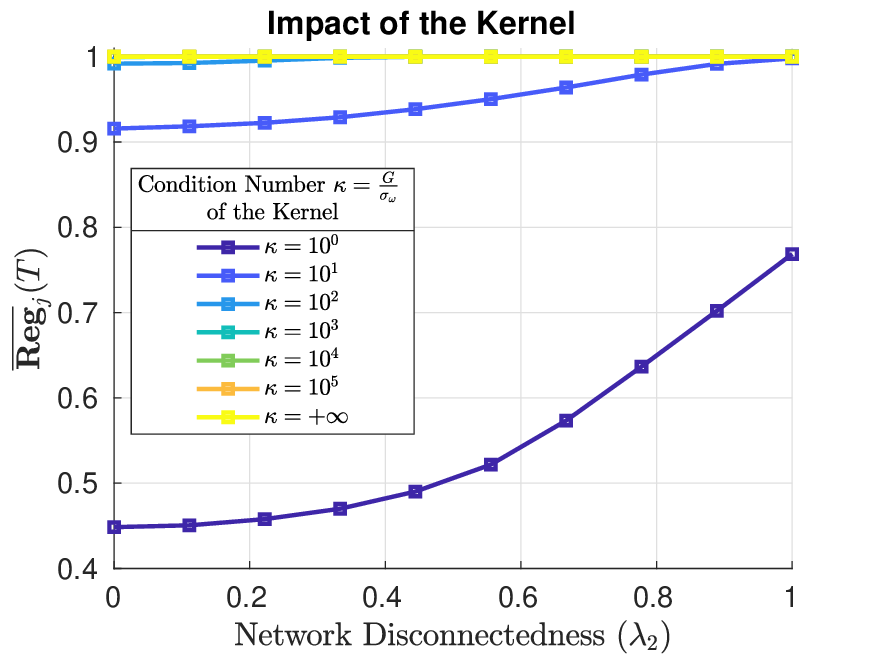}
    \caption{N-ISR against the network topology for various conditions numbers for the Kernel (with $T=5$).}
    \label{fig:conditions_numbers_kernel}
\end{figure}

DOMD and DAOL show similar PEP bounds w.r.t. $N$, while~\Cref{eq:bound_DOMD} is in $\mathcal{O}(n^2)$.
In Fig.~\ref{fig:conditions_numbers_kernel}, we observe that the condition number 
$\kappa = \frac{G}{\sigma_{\omega}}$ of the kernel substantially influences the extent 
to which agents benefit from communication within the network. For well-conditioned 
kernels, the regret is reduced by more than half, underscoring the importance of carefully 
selecting kernel functions in decentralized settings. Notably, this effect surpasses the 
improvements obtained by optimizing step-sizes, as discussed in~\Cref{section:improvement_and_design}.

\section{Improvement and design of methods}\label{section:improvement_and_design}

As a proof-of-concept, we now show how PEP allows for choosing better parameters for DOO methods i.e. to minimize the worst-case bound, by choosing step-sizes. For this purpose, we optimize DAOL and DOCG since they leave more latitude over the choice of step-sizes. Similar parameters optimization strategies can also be used on DOMD.

\subsection{Set-up}

\begin{figure}[htpb]
    \centering
    \includegraphics[width=0.75\linewidth]
    {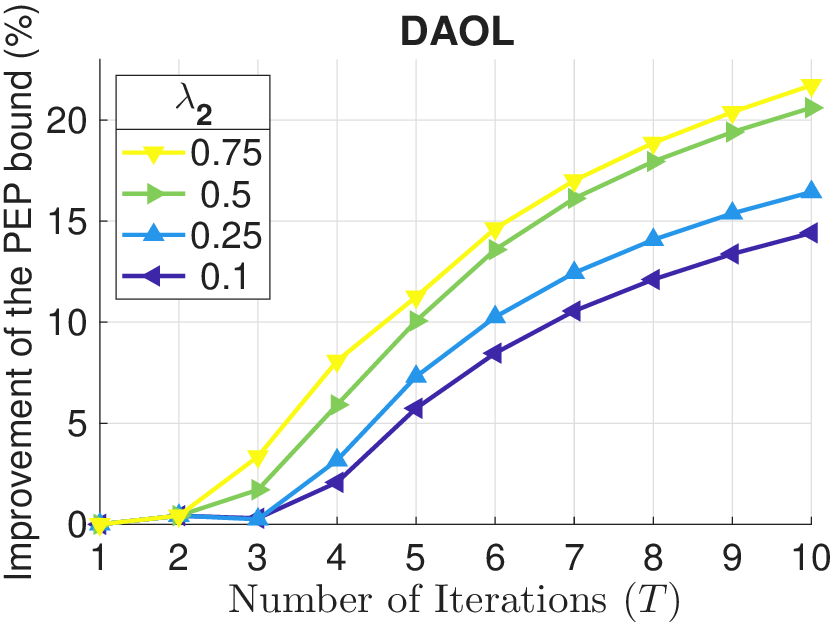}
    \vspace{0.5em}
    \includegraphics[width=0.75\linewidth]{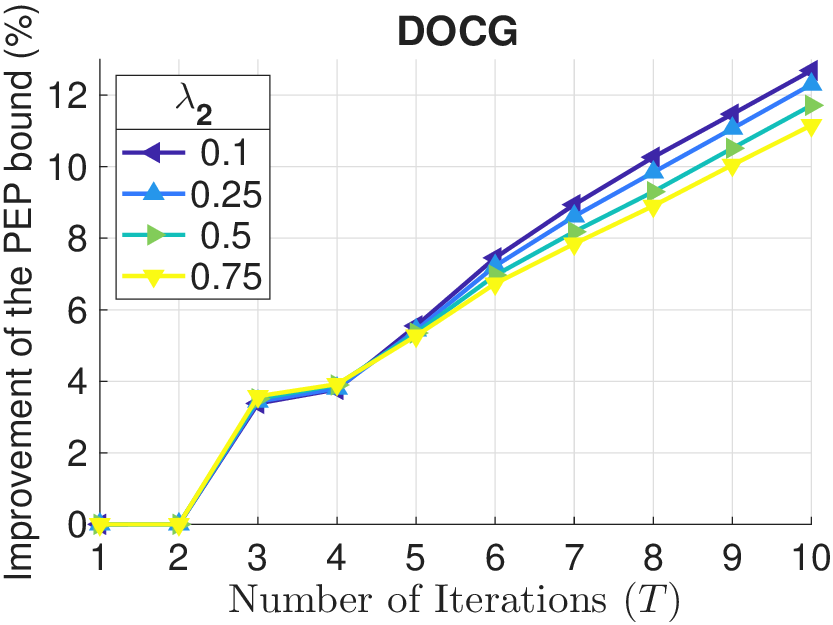}
    \caption{Improvement of DAOL and DOCG for various network classes.}
    \label{fig:design_method}
\end{figure}

We use the same parameters as in~\Cref{sub_section:set_up}.
Considering the following scheme for the step-sizes $\eta_t(a,b,c) := \frac{a}{t^c +b}$ which encompasses the step-sizes strategy from the literature for both DAOL and DOCG, we want to find parameters $a \in \left[0,D\right], \, b \in \left[0,L\right], \, c \in \left[0,1\right]$ so as to minimize the worst-case bound given by PEP. To do so, we use PEP as an oracle $\mathcal{P}_T^*(a,b,c) := \texttt{worst-case}(\eta_t(a,b,c)) / TnLD $ and $(a,b,c)\mapsto\mathcal{P}_T^*(a,b,c)$ is then minimized thanks to a black-box surrogate optimization algorithm (\texttt{surrogateopt} in MATLAB). 
The solver is warm-started for each $T \in \left[10\right]$ by the previous improved step-size strategy for $T-1$ iterations. For $T=1$, our initial guess if to use the step-sizes from the literature.
For each $T$, we allow $50$ PEP evaluations by \texttt{surrogateopt} and we call $(a_T^*,b_T^*,c_T^*)$ the triplet of improved parameters returned by the black-box solver after the number of allowed iterations. Our approach is tested on various network topologies $\lambda_2 \in \set{0.1,0.25,0.5,0.75}$ and results are shown in Fig.~\ref{fig:design_method} in which the following quantity is plotted on the $y$-axis 
\begin{equation*}
    \texttt{Improvement}(T) = 100\left(1 -\texttt{PEP\_Opt}(T)/\texttt{PEP\_Lit}(T)\right) 
\end{equation*}
where $\texttt{PEP\_Opt}(T):=\mathcal{P}_T^*(a_T^*,b_T^*,c_T^*)$ is the near tight N-ISR returned by PEP for optimized step-sizes while $\texttt{PEP\_Lit}(T)$ is the near tight N-ISR returned by PEP for steps-sizes coming from the literature.
\vspace{-0.5em}
\subsection{Outcomes}

For all topologies, DAOL with improved step-sizes consistently outperforms DAOL using the step-size strategy proposed in the literature.  In addition, for $T=10$, one can expect to save up to $22\%$ over the regret. Since the improvement is increasing with $T$, we suspect that the we could even do better for larger $T$ or by allowing $\eta_t$ to be independent from any constant. Interestingly, we can see that there is more room for improvement when the communication network is poorly connected (e.g. $\lambda_2=0.75$), which is again typical, than for well connected graphs (e.g. $\lambda_2=0.1$). 

Akin to DAOL, DOCG is successfully improved by up to $12 \% $ for $T=10$. Unlike DAOL, DOCG tends to see a slightly better improvement for well-connected networks.
Furthermore, if we consider the discussion in~\Cref{section:improvement_and_design}, DOCG equipped with improved step-sizes will then benefit from the communication network in the worst-case. Again, this call for further and much more extensive investigations.

\section{CONCLUSION AND FUTURE WORK}

We have shown how the PEP framework allows to get near-tight worst-case bounds for DOO algorithms. Then, we have demonstrated that picking the ``best-method'' by using analytical bounds, which turn out to be conservative, can lead to choosing a sub-optimal method.
In fact, in the chosen range of parameters, analytical bounds would have suggested to use DOCG instead of DAOL and DOMD. It turns out that PEP bounds tend to show it would be more reasonable to use DAOL or DOMD over DOCG.
In addition, surprising results concerning the performance of DOCG w.r.t. the network topology were observed and call for future investigations.

Furthermore, DAOL and DOCG have been improved by choosing the step-sizes so as to minimize the near-tight worst-case Regret. Our approach can be repeated on other DOO methods. In our case, it is worth noticing that tight bounds have been improved by up to about $20\%$.

We emphasize that these results are numerical and were obtained for a specific range / choice of parameters. They suggest the potential for very significant improvements in the general analysis. As future work, an in-depth study must be performed with wider ranges of parameters as well as deriving analytical formulas for the worst-case bounds.

\addtolength{\textheight}{-12cm} 

\vspace{-0.5em}
\section*{APPENDIX}\label{Appendix}

PEP models for DOO algorithms can be found in {\url{github.com/ErwanMeunier/PEP-for-Decentralized-Online-Optimization.git}}.
\vspace{-1em}
\AtNextBibliography{\small}
\printbibliography

\end{document}